\documentclass[11pt,bezier]{article}
\usepackage{amsmath,amssymb,amsfonts,euscript,graphicx}

\title{\bf   Prime {\it M}-Ideals, {\it M}-Prime Submodules,    {\it M}-Prime Radical   and  {\it M}-Baer's  Lower Nilradical
 of Modules \thanks {The research
 of the second author was in part supported by
a grant from IPM (No. 90160034).}
\thanks {2010{ \it Mathematics Subject Classification}: 16S38, 16D50, 16D60,
  16N60. }}

\author{{ {\bf J. A.  Beachy}$^{{\rm a}\flat}$,  {\bf M. Behboodi}$^{{\rm b,c}\natural}$  {\bf
and} {\bf  F. Yazdi}$^{{\rm b}\sharp}$}\\
 {\small{ $^{\rm a}$Department of Mathematical
 Sciences,  Northern Illinois
 University}}\vspace{-1mm}\\ {\small{ DeKalb, IL, USA 60115-2888}}\\
 {\small{ $^{\rm b}$Department of Mathematical Sciences,  Isfahan
University of Technology}}\vspace{-1mm}\\ {\small{ Isfahan,
Iran, 84156-83111}}\\
{\small{ $^{\rm c}$School of Mathematics, Institute for Research
in
Fundamental Sciences (IPM)}}\vspace{-1mm}\\ {\small{  Tehran, Iran, 19395-5746}}\vspace{-1mm}\\
 {\footnotesize{$\mathsf{{^\flat}beachy@math.niu.edu}$}}\vspace{-1mm}\\
{\footnotesize{$\mathsf{^{\natural}mbehbood@cc.iut.ac.ir}$}}\vspace{-1mm}\\
{\footnotesize{$\mathsf{^{\sharp}f.yazdi@math.iut.ac.ir}$}}}
\date{}
\begin{document}
\maketitle
\begin{abstract}
{\small \noindent Let $M$ be a fixed left $R$-module. For a left
$R$-module $X$, we introduce the notion of
  $M$-prime (resp. $M$-semiprime) submodule   of $X$ such that in the case $M=R$, which coincides with
   prime (resp. semiprime) submodule of $X$. Other concepts encountered in the general theory
   are $M$-$m$-system sets, $M$-$n$-system sets, $M$-prime radical and M-Baer's lower nilradical of  modules. Relationships
between these concepts and basic properties are established. In
particular, we identify certain submodules of $M$, called ``prime
$M$-ideals'', that play a role analogous to that of prime
(two-sided) ideals in the ring $R$. Using this definition, we show
that if $M$ satisfes condition $H$ (defined latter) and ${\rm
Hom}_R(M,X)\neq 0$ for all modules $X$ in the category
$\sigma[M]$, then there is a one-to-one correspondence between
isomorphism classes of indecomposable $M$-injective modules in
$\sigma[M]$ and prime $M$-ideals of $M$.   Also, we investigate
the prime $M$-ideals, $M$-prime submodules and $M$-prime radical
of Artinian modules.}
\end{abstract}

\section{Introduction}

 All rings in this paper are associative with identity and modules are unitary left modules.
   Let $R$ be a ring and $X$ be an $R$-module. If $Y$ is a submodule (resp. proper submodule)
    of $X$ we write $Y\leq X$ (resp. $Y\lneqq X$).

In the literature, there are many different generalizations of the
notion of prime two-sided  ideals to  left ideals  and also to
modules. For instance, a   proper left ideal $L$ of a ring $R$ is
called prime if, for any  elements $a$ and $b$ in $R$ such that
$aRb\subseteq L$, either $a\in L$ or $b\in L$.   Prime left ideals
have properties
 reminiscent of prime ideals in commutative rings. For example,
 Michler \cite{Mic} and Koh \cite{Koh} proved that the ring $R$ is left
 Noetherian if and only if every prime left ideal is finitely
 generated. Moreover, Smith \cite{Smi},  showed that if $R$ is left
 Noetherian (or even if $R$ has  finite left Krull dimension) then a left
 $R$-module $X$ is injective if and only if, for every essential
 prime left ideal $L$ of $R$ and homomorphism $\varphi :
 L\rightarrow X$, there exists a homomorphism $\theta :
 R\rightarrow X$ such that $\theta|_L=\varphi$.
 Let us mention another generalization of the notion of prime
 ideals to  modules.
 Let $X$ be a left $R$-module. If $X\neq 0$ and
${\rm Ann}_R(X)={\rm Ann}_R(Y)$ for all nonzero submodules $Y$ of
$X$ then $X$ is called  a {\it prime module}. A proper submodule
$P$ of $X$ is called a {\it prime submodule} if  $X/P$ is a prime
module, i.e., for every ideal $I\subseteq R$ and every submodule
$Y\subseteq X$, if $IY\subseteq P$, then either $Y\subseteq P$ or
$IX\subseteq P$. The notion of prime submodule was first
introduced and systematically studied by Dauns \cite{Dau} and
recently has received some attention. Several authors have
extended the theory of prime ideals of $R$ to prime submodules,
(see \cite{Beh1,Beh2,Beh3,Dau,JS,Lu,MM2,MS} ). For example, the
classical result of Cohen's is extended to prime submodules over
commutative rings, namely a finitely generated module is
Noetherian if and only if every prime submodule is finitely
generated (see \cite[Theorem 8]{Lu} and \cite{Kar}) and also any
Noetherian module contains only finitely many minimal prime
submodules (see \cite[Theorem 4.2]{MS}).

We assume throughout the paper $_RM$ is a fixed left $R$-module.
The category $\sigma[M]$ is defined to be the full subcategory of
$R$-Mod that contains all modules $_RX$ such that $X$ is
isomorphic to a submodule of an $M$-generated module (see
\cite{Wis} for more detail).

 Let ${\cal {C}}$ be  a class of modules in $R$-Mod, and let $\Omega$
 be the set of kernels of $R$-homomorphisms from $M$ in to ${\cal{C}}$. That
 is, $$\Omega=\{K \subseteq M ~\mid ~\exists ~W \in{\cal {C}}~and
~f\in {\rm Hom}_R(M, W)~ with~ K = {\rm ker} (f)\}.$$ Then the
{\it annihilator of ${\cal {C}}$ in $M$}, denoted by ${\rm
Ann}_M({\cal{C}})$, is defined to be the intersection of all
elements of $\Omega$, i.e.,
              ${\rm Ann}_M({\cal{C}})=\bigcap_{K\in\Omega}K$.

  Let $N$ be a submodule of $M$. Following Beachy \cite{Beachy}, for
each module $_RX$ we define  $$N\cdot X={\rm Ann}_X({\cal {C}}),$$
where ${\cal{C}}$ is the class of modules $_RW$ such that
$f(N)=(0)$ for all $f\in {\rm Hom}_R(M, W)$. It follows
immediately from the definition that
\begin{center}
     $N\cdot X=(0)$ if and only if $f(N) = (0)$ for all $f\in {\rm Hom}_R(M, X)$.
     \end{center}
~~~~Clearly the class ${\cal{C}}$ in definition of $N\cdot X$ is
closed under formation of submodules and direct products, and so
$N\cdot X$ is the smallest submodule $Y \subseteq X$ such that
$N\cdot (X/Y) = (0)$.

The submodule $N$ of $M$ is called an {\it $M$-ideal} if there is
a class ${\cal{C}}$ of modules in $\sigma[M]$ such that $N={\rm
Ann}_M({\cal{C}})$. Note that although the definition of an
$M$-ideal is given relative to the subcategory $\sigma[M]$, it is
easy to check that $N$ is an $M$-ideal if and only if $N={\rm
Ann}_M({\cal{C}})$ for some class ${\cal{C}}$  in $R$-Mod (see
\cite[Page 4651]{Beachy}).

 In this article for a left  $R$-module $X$, we introduce the notions of $M$-prime
 submodule,  $M$-semiprime
 submodule of $X$ and prime $M$-ideal of $M$ as follows:

\noindent{\bf Definition 1.1.} Let $X$ be an  $R$-module. A proper
submodule  $P$ of   $X$ is called an {\it $M$-prime submodule} if
for all submodules  $N\leq M,~ Y\leq
 X$, if $N\cdot Y\subseteq P$,    then
  either   $N\cdot X \subseteq P$ or $Y \subseteq P$. An $R$-module $X$ is called an {\it $M$-prime
  module}   if  $(0)\lneqq X$ is an  $M$-prime submodule.  Also, a proper submodule
 $P$ of $X$ is called an {\it $M$-semiprime submodule} if  for all submodules  $N\leq M,~ Y\leq
 X$, if $N^2\cdot Y\subseteq P$,    then $N\cdot Y\subseteq P$, where
 $N^2:=N\cdot N$. An $R$-module $X$ is called an {\it $M$-semiprime
 module}     if  $(0)\lneqq X$ is an  $M$-semiprime submodule.

\noindent{\bf Definition 1.2.}  A proper $M$-ideal $P$  of $M$ is
called a  {\it prime $M$-ideal}  (resp. {\it semiprime $M$-ideal})
if there exists an $M$-prime module  (resp. $M$-semiprime module)
$_RX$ such that $P={\rm Ann}_M(X)$.

It is clear that in case $M=R$, the notion of an $R$-prime
submodule (resp. $R$-semiprime submodule) reduces to the familiar
definition of a prime submodule (resp. semiprime submodule). Also,
the notion of an $R$-ideal (resp.  prime $R$-ideal) of $_RR$
reduces to the familiar definition of an ideal (resp. a  prime
ideal) of $R$.

Recently, the idea of {\it $M$-prime module}  was introduced and
extensively studied by Beachy [1] by defining a module $_RX$ to be
{\it $M$-prime} if ${\rm Hom}_R(M, X)\neq 0$, and ${\rm
Ann}_M(Y)={\rm Ann}_M(X)$ for all submodules $Y\subseteq X$ such
that ${\rm Hom}_R(M, Y)\neq 0$. Also, he defined an
  $M$-ideal $P$ to be {\it prime $M$-ideal} if there exists an
$M$-prime module $_RX$ such that $P={\rm Ann}_M(X)$. Clearly, our
definition of $M$-prime module is slightly different than Beachy,
and hence, for the sake of clarity, for the remainder of the paper
we will use the term ``Beachy-$M$-prime module" (resp.
``Beachy-prime $M$-ideal") rather than ``$M$-prime module" (resp.
``prime $M$-ideal") of Beachy \cite{Beachy}, respectively.

 In ring theory,
prime ideals are closely tied  to m-system  sets (a nonempty set
$S \subseteq R$ is  said to be an {\it $m$-system set} if for each
pair
 $a, b$ in  $S$, there exists $r\in  R$ such that $arb \in S$). The
complement of a prime ideal is an  m-system, and given an m-system
set $S$, an ideal disjoint from $S$ and maximal with respect to
this property is always a prime ideal. Moreover, for an ideal $I$
in a ring $R$, the set  $\sqrt{I}:=\{ s\in R \ |$  every m-system
containing $s$ meets $I\}$ equals the intersection of all the
prime ideals containing $I$. In particular, $\sqrt{I}$ is a
semiprime ideal in $R$ and $\sqrt{(0)}$ is called {\it Baer-McCoy
radical} (or {\it prime radical}) of $R$ (see for example
\cite[Chapter 4]{Lam}, for more details). In this paper, we extend
these facts for $M$-prime submodules. Relationships between these
concepts and basic properties are established. In Section 2, among
other results, for an $R$-module $X$ we define {\it $M$-Baer-McCoy
radical} (or {\it $M$-prime radical}) of $X$, denoted ${\rm
rad}_M(X)=\sqrt[M]{(0)}$,  to be the intersection of of all the
$M$-prime submodules in  $X$. Also, in Section 3,  we extend the
notion of nilpotent and strongly nilpotent element of modules to
$M$-nilpotent and strongly $M$-nilpotent element  of modules $X\in
\sigma[M]$ for a fix module $M$. Also,   for an $R$-module $X\in
\sigma[M]$, we define {\it $M$-Baer's lower nilradical} of $X$,
denoted by $M$-${\rm Nil}_*(_RX)$, to be the set of all strongly
$M$-nilpotent elements of $X$. In particular, it is shown that  if
$M$ is  projective  in $\sigma[M]$, then for each $X\in
\sigma[M]$, ${\rm Nil}_*(M).X\subseteq M$-${\rm
Nil}_*(_RX)\subseteq {\rm rad}_M(X)$ (see Proposition 3.6).

In Section 4, we rely on the prime $M$-ideals of $M$ that play a
role analogous to that of prime ideals in the ring $R$. The module
$_RX$ is called {\it $M$-injective} if each $R$-homomorphism $f:
K\rightarrow X$ defined on a submodule $K$ of $M$ can be extended
to an $R$-homomorphism $\widehat{f}: M\rightarrow  X$ with $f=
\widehat{f}i$, where $i: K \rightarrow M$ is the natural inclusion
mapping. We note that Baer's criterion for injectivity shows that
any $R$-injective module is injective in the category $R$-Mod of
all left $R$-modules. It is well-known that if $R$ is a
commutative Noetherian ring, then there is a one-to-one
correspondence between isomorphism classes of indecomposable
injective $R$-modules and prime ideals of $R$. Gabriel showed in
\cite{Gab} that this one-to-one correspondence remains valid for
any left Noetherian ring that satisfies what he called ''condition
$H$''. In current terminology, a module $_RX$ is said to be
finitely annihilated if there is a finite subset $x_1,\cdots,x_n$
of $X$ with ${\rm Ann}_R(X)={\rm Ann}_R(x_1,\cdots,x_n)$. Then by
definition the ring $R$ satisfies condition $H$ if and only if
every cyclic left $R$-module is finitely annihilated. It follows
immediately that, the ring $R$ satisfies condition $H$ if and only
if every finitely generated left $R$-module is finitely
annihilated. We note the stronger result due to Krause \cite{Krau}
that if $R$ is left Noetherian, then there is a one-to-one
correspondence between isomorphism classes of indecomposable
injective left $R$-modules and prime ideals of $R$ if and only if
$R$ is a left fully bounded ring (see  \cite[Theorem 8.12]{Good}
for a proof). In \cite[Theorem 6.7]{Beachy},  Beachy shown that
Gabriel's correspondence can be extended to $M$-injective modules,
provided that ${\rm Hom}_R(M, X)\neq 0$ for all modules $X$ in
$\sigma[M]$. In Section 4, by using our definition of prime
$M$-ideal, we show that also there is Gabriel's correspondence
between indecomposable $M$-injective modules in $\sigma[M]$ and
our prime $M$-ideals.

Finally, in Section 5, we study the prime $M$-ideal, $M$-prime
submodules and $M$-prime radical of Artinian modules. The {\it
prime radical} of the module $M$, denoted by $P(M)$, is defined to
be the intersection of all prime $M$-ideals of $M$. Recall that a
proper submodule $P$ of $M$ is {\it virtually maximal} if the
factor module $M/P$ is a homogeneous semisimple $R$-module, i.e.,
$M/P$ is a direct sum of isomorphic simple modules.  It is shown
that if $M$ is  an Artinian  $M$-prime module, then $M$ is a
homogeneous semisimple module (see Proposition 5.1). In
particular, if $M$ is an Artinian $R$-module such that it is
projective  in $\sigma[M]$, then every prime $M$-ideal of $M$ is
virtually maximal and $M/P(M)$ is a Noetherian $R$-module (see
Theorem 5.6). Moreover, either $P(M)=M$ or there exist primitive
(prime) $M$-ideals $P_1$,...,$P_n$ of $M$ such that
$P(M)=\bigcap_{i=1}^{n}P_i$ (see Theorem 5.7).

\section{$M$-prime submodules and $M$-prime radical of modules }

We begin this section with the following three useful lemmas.

 \noindent{\bf
Lemma 2.1.} [1, Proposition 1.6] {\it Let $N$ be a submodule of
$M$. Then for any $R$-module $X$, $N\cdot X=(0)$ if and only if
$N\subseteq Ann_M(X)$.}

\noindent{\bf Lemma  2.2.} [1,  Proposition 1.9] {\it Let $N$ and
$K$ be submodules of $M$.}\vspace{3mm}\\
(a) {\it If $N\subseteq K$, then $N\cdot X\subseteq K\cdot X$ for
all submodules $_RX$.}\vspace{1mm}\\
(b) {\it If $K$ is an $M$-ideal, then so is  $N\cdot K$.}\vspace{1mm}\\
(c) {\it The submodule  $N\cdot M$ is the smallest $M$-ideal that
contains $N$.}\vspace{1mm}\\
(a) {\it If $N$ is an $M$-ideal, then $N\cdot K\subseteq N\cap
K$.}

\noindent{\bf Lemma  2.3.}  {\it Let $Y_1$, $Y_2$ be submodules of
$_RX$. If $Y_1\subseteq Y_2$, then $N\cdot Y_1\subseteq N\cdot
Y_2$, for each  submodule $N$ of $M$.}

\noindent {\bf Proof.}  Suppose $N\leq M$ and $Y_1$, $Y_2$ be
submodules of $_RX$ with $Y_1\subseteq Y_2$. Then  $N\cdot
Y_1={\rm Ann}_{Y_1}({\cal {C}})$  and  $N\cdot Y_2={\rm
Ann}_{Y_2}({\cal {C}})$, where ${\cal{C}}$ is the class of modules
$_RW$ such that $f(N)=(0)$ for all $f\in {\rm Hom}_R(M, W)$. On
the other hand  $N\cdot Y_i=\bigcap _{K\in\Omega_i}K$  $(i=1,2)$,
where
$$\Omega_i=\{K \subseteq Y_i ~\mid ~\exists ~W \in{\cal {C}}~and
~f\in {\rm Hom}_R(Y_i, W)~ with~ K = {\rm ker} (f)\}$$ Clearly,
for each $f\in {\rm Hom}_R(Y_2, W)$, $f|_{Y_1}\in {\rm Hom}_R(Y_1,
W)$, where $f|_{Y_1}$ is the restriction of $f$ on $Y_1$. Since
${\rm ker}(f|_{Y_1})\subseteq {\rm ker}(f)$, we conclude that for
each $K\in \Omega_2$,  there exists $K'\in \Omega_1$ such that
$K'\subseteq K$. Thus $N\cdot Y_1\subseteq N\cdot Y_2$.  $\Box$

The following evident proposition offers several characterizations
of an  $M$-prime module.

\noindent{\bf Proposition 2.4.}  {\it Let $X$ be a nonzero
$R$-module.
Then the following statements are equivalent.}\vspace{3mm}\\
(1) {\it $X$ is an $M$-prime module.}\vspace{1mm}\\
(2) {\it For every submodule $N\subseteq M$ and every nonzero
submodule $Y\subseteq X$, if $N\cdot Y=(0)$, \indent then   $N\cdot X=(0)$.}\vspace{1mm}\\
(3) {\it For every $M$-ideal $N\subseteq M$ and every nonzero
submodule $Y\subseteq X$, if $N\cdot Y=(0)$, then \indent  $N\cdot X=(0)$.}\vspace{1mm}\\
(4) {\it For every nonzero submodules  $Y_1$, $Y_2\subseteq X$,
$Ann_M(Y_1)=Ann_M(Y_2)$.}\vspace{1mm}\\
(5) {\it  Every nonzero submodule  $Y\subseteq X$ is an $M$-prime
module.}\vspace{1mm}\\
(6) {\it  For every nonzero  submodule  $Y\subseteq X$,
$P=Ann_M(Y)$ is a prime $M$-ideal of $M$  and \indent
$P=Ann_M(X)$.}

\noindent {\bf Proof.}   $(1)\Rightarrow (2)\Rightarrow (3)$ is  clear.\\
$(3)\Rightarrow (4)$. Let $Y_1 , Y_2$ be two nonzero submodules of
$X$ and let $N_1 := {\rm Ann}_M(Y_1)$,  $N_2 := {\rm Ann}_M(Y_2)$.
Thus by Lemma 2.1,  $N_1\cdot Y_1=(0)$ and $N_2\cdot Y_2=(0)$.
Since $N_1$, $N_2$ are $M$-ideals, $N_1\cdot X=N_2\cdot X=(0)$ by
$(3)$. Thus  $N_1\subseteq {\rm Ann}_M(X)$ and $N_2\subseteq {\rm
Ann}_M(X)$. On the other hand ${\rm Ann}_M(X)\subseteq N_1$ and
${\rm Ann}_M(X)\subseteq N_2$.
Thus $N_1=N_2={\rm Ann}_M(X)$.\\
$(4)\Rightarrow (5)$. Let $Y$ be a nonzero submodule of $X$.
Assume that $N$ is  a submodule of $M$ and $Z$ be a nonzero
submodule of Y such that $N\cdot Z=(0)$. So $N\subseteq {\rm
Ann}_M(Z)$. By (4),  ${\rm Ann}_M(Z)={\rm Ann}_M(X)$  and so it
follows that $N\subseteq {\rm Ann}_M(X)$ and hence $N\cdot X=(0)$.
Since $N\cdot Y\subseteq N\cdot X$,  so $N\cdot Y=(0)$. Thus $Y$ is an $M$-prime module.\\
$(5)\Rightarrow (1)$ and  $(5)\Rightarrow (6)\Rightarrow(4)$ are
clear. $\Box$

\noindent{\bf Remark  2.5.} Clearly every simple $R$-module $X$ is
an $M$-prime module. Now let $R$ be a domain which is not a field
and let $M$ be a nonzero divisible $R$-module. Then every nonzero
simple $R$-module $X$ is an $M$-prime module, but  $X$ is not a
Beachy-$M$-prime module, since ${\rm Hom}_R(M, X)=0$.

   The following lemma shows that in the case  ${\rm Hom}_R(M,X)\neq 0$, if  $X$ is an
   $M$-prime module then  $X$ is also a Beachy-$M$-prime module.

\noindent{\bf Lemma  2.6.} ([1,  Proposition 2.2])  {\it Let $X$
be an  $R$-module such that ${\rm Hom}_R(M,X)\neq 0$.
Then the following statements are equivalent.}\vspace{3mm}\\
(1) {\it $X$ is a Beachy-$M$-prime module.}\vspace{1mm}\\
(2) {\it For every $M$-ideal $N$ of $M$ and every nonzero
submodule $Y$ of $X$ with $M\cdot Y\neq (0)$,  if  \indent $N\cdot Y=(0)$, then   $N\cdot X=(0)$.}\vspace{1mm}\\
(3) {\it For each $m\in M\setminus Ann_M(X)$ and  each $0\neq f\in
Hom_R(M, X)$,    there exists   $g\in \indent Hom_R(M, f(M))$ such
that
$g(m)\neq 0$.}\vspace{1mm}\\
(4) {\it  For any $M$-ideal $N\subseteq M$ and any $M$-generated
submodule $Y\subseteq X$, if  $N\cdot Y=(0)$, then  \indent
$N\cdot X=(0)$.}

\noindent{\bf Proposition 2.7.} {\it Let $X$ be an $R$-module such
that $Hom_R(M,X)\neq 0$. If  $X$ is an $M$-prime module then  $X$ is a
 Beachy-$M$-prime module.}

 \noindent {\bf Proof.}  By Proposition 2.4 and Lemma 2.6, is clear. $\Box$

 The following example shows that the   converse of
Proposition 2.7 is not true in general.

   \noindent{\bf Example  2.8.} Let $R=\Bbb{Z}$. For each prime number $p$, ${\rm Hom}_{\Bbb{Z}}({\Bbb{Z}}_{p_{\infty}}, {\Bbb{Z}}_{p_{\infty}})\neq 0$ and
   for each proper $\Bbb{Z}$-submodule $Y\subsetneqq {\Bbb{Z}}_{p_{\infty}}$, ${\Bbb{Z}}_{p_{\infty}}\cdot Y=(0)$, since
    ${\rm Hom}_{\Bbb{Z}}({\Bbb{Z}}_{p_{\infty}}, Y)=(0)$.
   Thus  by Lemma 2.6,  ${\Bbb{Z}}_{p_{\infty}}$ is a Beachy-${\Bbb{Z}}_{p_{\infty}}$-prime module but it is not  a ${\Bbb{Z}}_{p_{\infty}}$-prime module,  since
   ${\Bbb{Z}}_{p_{\infty}}\cdot {\Bbb{Z}}_{p_{\infty}}\neq(0)$.

\noindent{\bf Lemma  2.9.} ([1,  Proposition 5.5]) {\it Assume
that $M$ is projective in $\sigma[M]$, and let $N$ be any
submodule of $M$. The following conditions hold for any module
$_RX$ in $\sigma[M]$ and any submodule $Y\subseteq X$.}\vspace{3mm}\\
(a) {\it $N\cdot X =\sum_{f\in Hom_R(M,X)}f(N)$.}\vspace{1mm}\\
(b) {\it $N\cdot (X/Y ) =(0)$ if and only if $N\cdot X\subseteq Y$.}\vspace{1mm}\\
(c) {\it If $N=Ann_M(X/Y)$, then $Ann_M(X/(N\cdot X))= N$.}

\noindent{\bf Proposition 2.10.} {\it Assume that $M$ is
projective in $\sigma[M]$,  and let $_RX \in \sigma[M]$. Then:}\vspace{3mm}\\
(i) {\it  For a submodule $P \lneqq X$,  if $P$ is an  $M$-prime
submodule of $X$, then
 $X/P$ is  an $M$-\indent prime module.}\vspace{2mm}\\
(ii) {\it For an $M$-ideal  $P \lneqq M$,  the following
conditions are equivalent.}\vspace{3mm}\\
\indent (1) {\it $P$ is a prime $M$-ideal.}\vspace{1mm}\\
\indent (2) {\it $P$ is an  $M$-prime
submodule of $M$.}\vspace{1mm}\\
\indent (3) {\it $M/P$ is an $M$-prime module.}

 \noindent {\bf Proof.}  (i).  Let $N$ be a submodule of $M$
and $Y/P$ be a nonzero submodule of $X/P$ such that $N\cdot
(Y/P)=(0)$. By Lemma 2.9 (b),  $N\cdot Y\subseteq P$. Since $P$ is
an $M$-prime submodule,  either $N\cdot X\subseteq P$ or
$Y\subseteq P$. If $Y\subseteq P$, then $Y/P=(0)$,  a
contradiction. Thus $N\cdot X\subseteq P$ and so $N\cdot
(X/P)=(0)$ by Lemma 2.9 (b). Thus  by Proposition 2.4, $X/P$ is an
$M$-prime module.\\
(ii) $(1)\Rightarrow (2)$. Suppose that $P$ is a prime $M$-ideal
and $N\cdot K\subseteq P$, for an $M$-ideal $N$ and submodule $K$
of $M$ with $K\nsubseteq P$. By assumption there is an $M$-prime
module $X$ with $P={\rm Ann}_M(X)$, and so there exists $f\in {\rm
Hom}_R(M/P,X)$ with $f((K + P)/P)\neq (0)$. Since $N\cdot
K\subseteq P$, we have $N\cdot K \subseteq P\cap K$. Now Lemma 2.9
(b) implies that $N\cdot (K/(P\cap K))=(0)$ and hence $N\cdot f((K
+ P)/P)=(0)$ (since $(K + P)/P\cong K/(P\cap K)$). Since $X$ is an
$M$-prime module,   $N\cdot X = (0)$ by Proposition 2.4, and so
$N\subseteq P$ (since
$P={\rm Ann}_M(X)$).\\
$(2)\Rightarrow (3)$. Let $N$ be an $M$-ideal and $K/P$ be a
nonzero submodule of $M/P$ such that $N\cdot (K/P)=(0)$. Since $M$
is projective in $\sigma[M]$, so $N\cdot K\subseteq P$ by Lemma
2.9 (b). Now by (2) either  $N\subseteq P$ or $K\subseteq P$.
Since $K/P\neq (0)$, so $K\nsubseteq P$ and hence $N\subseteq P$.
On the other hand $N\cdot M=N$, since $N$ is an $M$-ideal. Thus
$N\cdot M\subseteq P$ and hence by Lemma 2.9 (b),  $N\cdot
(M/P)=(0)$. Now $M/P$ is an $M$-prime module by Proposition
2.4.\\
$(3)\Rightarrow (1)$. Since $P$ is an $M$-ideal,  $P={\rm
Ann}_M(M/P)$ and since  $M/P$ is  an $M$-prime module, we conclude
that  $P$ is a prime $M$-ideal. $\Box$

The following example shows that even in the case the $R$-module
$M$ is projective in $\sigma[M]$, an $M$-prime module need not be
a Beachy-$M$-prime module.

   \noindent{\bf Example  2.11.} Let $R=\Bbb{Z}$ and $M=\Bbb{Q}$ as $\Bbb{Z}$-module. Then it is
    easy to check that
    $\Bbb{Q}$ is projective in $\sigma[{\Bbb{Q}}]$. Clearly, for each prime number $p$,
    ${\Bbb{Z}}_p$ is a $\Bbb{Q}$-prime module, but it is not a Beachy-$\Bbb{Q}$-prime module,
    since ${\rm Hom}_{\Bbb{Z}}({\Bbb{Q}},
    {\Bbb{Z}}_p)=(0)$.

Now we have to adapt the notion of an $M$-m-system set  to modules
$_RX$ (Behboodi in \cite{Beh1}, has generalized the notion of
m-system of rings to modules).

\noindent{\bf Definition 2.12.} Let $X$ be  an $R$-module. A
nonempty set $S\subseteq X\setminus \{0\}$ is called an
$M$-$m$-system if, for each submodule $N\subseteq M$, and for all
submodules $Y,Z \subseteq X$, if $(Y+Z)\cap S\neq\emptyset$ and
$(Y + N\cdot X) \cap S \neq \emptyset$, then $(Y + N\cdot Z) \cap
S \neq \emptyset$.

\noindent{\bf Corollary 2.13.} {\it Let  $X$ be an $R$-module.
Then a submodule $P\lneqq X$ is $M$-prime if and only if
$X\setminus P$ is an $M$-$m$-system.}

\noindent {\bf Proof.}   $(\Rightarrow)$. Suppose $S=X\setminus
P$. Let $N$ be a submodule of $M$ and $Y$, $Z$ be submodules of
$X$ such that $(Y+Z)\cap S\neq \emptyset$ and $(Y+ N\cdot X)\cap
S\neq \emptyset$. If $(Y+ N\cdot Z)\cap S= \emptyset$ then $Y+
N\cdot Z\subseteq P$. Hence $N\cdot Z\subseteq P$ and since $P$ is
an $M$-prime submodule, $Z\subseteq P$ or $N\cdot X\subseteq P$.
It follows that $(Y+Z)\cap S=\emptyset$ or $(Y+ N\cdot X)\cap
S=\emptyset$, a contradiction. Therefore, $S \subseteq X \setminus \{0\}$ is an $M$-m-system set.\\
$(\Leftarrow)$. Let $S=X\setminus P$ be an $M$-m-system in $X$.
Suppose $N\cdot Z\subseteq P$, where $N$ is a submodule of $M$ and
$Z$ is a submodule $X$. If $Z\not\subseteq P$ and $N\cdot
X\not\subseteq P$, then $Z\cap S\neq\emptyset$ and $(N\cdot X)\cap
S\neq \emptyset$. Thus $(N\cdot Z)\cap S\neq \emptyset$, a
contradiction. Therefore, $P$ is an $M$-prime submodule of $X$.
$\Box$

\noindent{\bf Proposition 2.14.} {\it Let $X$ be an $R$-module,
$P$ be a proper submodule of $X$ and $S:= X\setminus P$. Then the
following statements are equivalent.}\vspace{3mm}\\
(1) {\it $P$ is  an $M$-prime submodule.}\vspace{1mm}\\
(2) {\it $S$ is an $M$-$m$-system.}\vspace{1mm}\\
(3) {\it For every  submodule $N\leq M$ and for every submodule $Z
    \leq X$, if $Z\cap S \neq \emptyset$ and \indent $(N\cdot X)\cap S \neq
    \emptyset$, then $(N\cdot Z)\cap S \neq \emptyset$.}

\noindent {\bf Proof.}   (1) $\Leftrightarrow$ (2) is by Corollary 2.13.\\
 (2) $\Rightarrow$ (3) is clear.\\
(3) $\Rightarrow$ (1). Suppose that $N\leq M$ and $Z\leq X$ such
that $N\cdot Z\subseteq P$. If $N\cdot X \nsubseteq P$ and $Z
\nsubseteq P$, then $(N\cdot X)\cap S \neq \emptyset$ and $Z\cap
S\neq\emptyset$. It follows that $(N\cdot Z)\cap S\neq\emptyset$
 by  (3), i.e.,  $N\cdot Z\nsubseteq P$, a contradiction.   $\Box$

\noindent{\bf Proposition 2.15.} {\it Let $X$ be an  $R$-module,
$S\subseteq X$ be an $M$-$m$-system  and $P$ be a submodule of $X$
maximal with respect to the property that $P$ is disjoint from
$S$. Then $P$ is an  $M$-prime submodule of $X$.}

\noindent{\bf  Proof.} Suppose $N\cdot Z\subseteq P$, where $N\leq
M$ and  $Z\leq X$. If $Z\not\subseteq P$ and $N\cdot
X\not\subseteq P$, then by the maximal property of $P$, we have,
$(P+Z)\cap S\neq\emptyset$ and $(P+ N\cdot X)\cap S\neq\emptyset$.
Thus $(P+ N\cdot Z)\cap S\neq\emptyset$ and it follows that $P\cap
S\neq\emptyset$, a contradiction. Thus $P$ must be an  $M$-prime
submodule.       $\Box$

Next we need a generalization of the notion of $\sqrt{Y}$ for any
submodule $Y$ of $X$. We adopt the
  following:

\noindent{\bf Definition 2.16.} Let $X$ be an  $R$-module. For a
submodule $Y$ of $X$, if there is an $M$-prime submodule
containing $Y$, then we define
\begin{center}
 $\sqrt[M]{Y}=\{x\in X:$ every $M$-m-system containing $x$ meets $Y\}.$
 \end{center}
  If there is no $M$-prime submodule containing $Y$,
then we put $\sqrt[M]{Y}=X$.

\noindent{\bf Theorem 2.17.} {\it  Let $X$ be an  $R$-module and
$Y\leq X$. Then either $\sqrt[M]{Y}=X$ or $\sqrt[M]{Y}$ equals the
intersection of all  $M$-prime submodules of $X$ containing $Y$.}

 \noindent {\bf Proof.}   Suppose that $\sqrt[M]{Y} \neq X$. This means that
\begin{center}
 $\{P :~P$  is an $M$-prime submodule of $X$ and $Y\subseteq P\}\neq\emptyset.$
\end{center}
We first prove that
\begin{center}
 $\sqrt[M]{Y}\subseteq\bigcap \{P :|~P$  is an $M$-prime submodule of $X$ and  $Y\subseteq P\}$.\\
\end{center}
Let $x \in \sqrt[M]{Y}$ and $P$ be any $M$-prime submodule of $X$
containing $Y$. Consider the $M$-$m$-system $X \setminus P$. This
$M$-$m$-system cannot contain $x$, for otherwise it meets $Y$ and
hence also $P$. Therefore, we have $x\in P$. Conversely, assume
$x\notin\sqrt[M]{Y}$. Then, by Definition 2.16, there exists an
$M$-m-system $S$ containing $x$ which is disjoint from $Y$. By
Zorn's Lemma, there exists a submodule $P\supseteq Y$ which is
maximal with respect to being disjoint from $S$. By Proposition
2.15, $P$ is an  $M$-prime submodule of $X$, and we have $x\notin P$,
 as desired.  $\Box$

Also, the following evident proposition offers several
characterizations   of $M$-semiprime modules.

\noindent{\bf Proposition 2.18.}  {\it Let $X$ be an  $R$-module.
Then the following statements are equivalent.}\vspace{3mm}\\
(1) {\it $X$ is an $M$-semiprime module.}\vspace{1mm}\\
(2) {\it For every submodule $N\subseteq M$ and every
submodule $Y\subseteq X$, if $N^2\cdot Y=(0)$,  then \indent  $N\cdot Y=(0)$.}\vspace{1mm}\\
(3) {\it Every nonzero submodule $Y\subseteq X$ is an
$M$-semiprime
module.}\vspace{1mm}\\
(4) {\it  For every nonzero  submodule  $Y\subseteq X$,
$P=Ann_M(Y)$ is a semiprime $M$-ideal.}

\noindent {\bf Proof.}   $(1)\Rightarrow (2)\Rightarrow (3)\Rightarrow (4)$ is clear.\\
$(4)\Rightarrow (1)$. Suppose $(0)\neq Y\leq X$ and  $N\leq M$
such that $N^2\cdot Y=(0)$. It follows that $N^2\subseteq {\rm
Ann}_M(Y)$ and since $P={\rm Ann}_M(Y)$ is a semiprime $M$-ideal,
there exists an $M$-semiprime module $Z$ such that ${\rm
Ann}_M(Y)={\rm Ann}_M(Z)$. Thus $N^2\cdot Z=(0)$ and so $N\cdot
Z=(0)$, i.e., $N\subseteq {\rm Ann}_M(Z)={\rm Ann}_M(Y)$. Thus
$N\cdot Y=(0)$. Therefore $X$ is an
 $M$-semiprime module.      $\Box$

\noindent{\bf Proposition 2.19.}  {\it Let $X$ be an  $R$-module.
Then any intersection of $M$-semiprime submodules of $X$ is an
$M$-semiprime submodule.}

\noindent {\bf Proof.}  Suppose  $Z_i\leq X$ $(i\in I)$ be
$M$-semiprime submodules of $X$ and put $Z=\bigcap_{i\in I}Z_i$.
Suppose $Y\leq X$ and $N\leq M$ such that $N^2\cdot Y\subseteq Z$.
It follows that $N^2\cdot Y\subseteq Z_i$ for each $i$. Since each
$Z_i$ is an $M$-semiprime submodule, $N\cdot Y\subseteq Z_i$ for
each $i$. Thus $N\cdot Y\subseteq Z$ and so $Z$ is an
$M$-semiprime submodule.  $\Box$

We recall the definition of the  notion of n-system in  a ring
$R$. A nonempty set $T \subseteq R$ is  said to be an {\it
$n$-system set} if for each  $a$ in  $T$, there exists $r\in  R$
such that $ara \in T$ (see for example \cite[Chapter 4]{Lam},  for
more details). The complement of a semiprime ideal is an  n-system
set, and if $T$ is an n-system in a ring $R$ such that $a\in T$,
then there exists an m-system $S\subseteq T$ such that $a\in S$
(see \cite[Lemma 10.10]{Lam}). This notion of n-system of rings
 has also generalized by Behboodi in \cite{Beh1} for modules. Now we have to adapt the notion of an
$M$-n-system set to modules $_RX$ .

\noindent{\bf Definition 2.20.} Let $X$ be an  $R$-module. A
nonempty set $T\subseteq X\setminus\{0\}$ is called an
$M$-$n$-system if, for every submodule $N\subseteq M$, and for all
submodules $Y,Z \subseteq X$, if $(Y+N\cdot Z)\cap
T\neq\emptyset$,
 then $(Y+N^2\cdot Z)\cap T\neq \emptyset$.

\noindent{\bf Proposition 2.21.} {\it Let $X$ be an $R$-module.
Then a submodule $P\lneqq X$ is an  $M$-semiprime submodule if and
only if $X\setminus P$ is an $M$-$n$-system.}

 \noindent {\bf Proof.}   $(\Rightarrow)$. Let $T=X\setminus P$. Suppose
 $N$ is a submodule of $M$ and $Y,Z$ are submodules of $X$ such
 that $(Y+N\cdot Z)\cap T\neq\emptyset$. If $(Y+N^2\cdot Z)\cap T=\emptyset$, then $(Y+N^2\cdot Z)\subseteq P$. Since $P$ is
 $M$-semiprime submodule, $(Y+N\cdot Z)\subseteq P$. Thus $(Y+N\cdot Z)\cap T=\emptyset$, a contradiction.
 Therefore, $T$ is an  $M$-$n$-system set in $X$.\\
 $(\Leftarrow)$. Suppose that $T=X\setminus P$ is an
 $M$-$n$-system in $X$. Suppose $N^2\cdot Z\subseteq P$, where $N\leq M$, $Z\leq X$,  but $N\cdot Z\nsubseteq P$.
  It follows that $(N\cdot Z)\cap T\neq\emptyset$ and so $(N^2\cdot Z)\cap T\neq\emptyset$, a contradiction.
  Therefore, $P$ is an $M$-semiprime  submodule of $X$.  $\Box$

The proof of the next proposition is similar to the proof of
Proposition 2.14.

\noindent{\bf Proposition 2.22.} {\it Assume that $P$ be a proper
submodule of $X$ and $T:=X\setminus P$. Then the
following statements are equivalent.}\vspace{3mm}\\
(1) {\it $P$ is an $M$-semiprime submodule.}\vspace{1mm}\\
(2) {\it  $T$ is an $M$-$n$-system set.}\vspace{1mm}\\
(3) {\it For every submodule $N \leq M$ and for every submodule
$Z\leq X$, if $(N\cdot Z)\cap T\neq\emptyset$, then \indent
 $(N^2\cdot Z)\cap T\neq\emptyset$.}

\noindent{\bf Lemma 2.23.} ([1, Proposition 5.6]) {\it Assume that
$M$ is projective in $\sigma[M]$, and let $K$, $N$ be submodules
of $M$. Then $(K\cdot N)\cdot X=K\cdot (N\cdot X)$ for any module
$_RX$ in $\sigma[M]$. }

\noindent{\bf Proposition 2.24.} {\it Assume that $M$ is
projective in $\sigma[M]$,  and let  $X\in \sigma[M]$. Then any
$M$-prime submodule of $X$ is an $M$-simiprime submodule.}

 \noindent {\bf Proof.}   Let $P \lneqq X$ be an $M$-prime
submodule of $X$ and $N\leq M$, $Y\leq X$ such that
$N^2.Y\subseteq P$. Since $M$ is projective in $\sigma[M]$, so
$N^2.Y=(N\cdot N).Y=N\cdot (N\cdot Y)$ by Lemma 2.23. Hence
$N\cdot (N\cdot Y)\subseteq P$. Now by assumption, $N\cdot
X\subseteq P$ or $N\cdot Y\subseteq P$. If $N\cdot Y\subseteq P$,
then $P$ is an $M$-semiprime submodule. If $N\cdot X\subseteq P$,
the  $N\cdot Y\subseteq N\cdot X\subseteq P$. Thus  $P$ is an
$M$-semiprime submodule.   $\Box$

\noindent{\bf Corollary 2.25.}  {\it Assume that $M$ is projective
in $\sigma[M]$ and $X\in \sigma[M]$. Then any intersection of
$M$-prime submodules of $X$ is an $M$-semiprime submodule.}

\noindent {\bf Proof.} Is by Proposition 2.19 and Proposition
2.24.  $\Box$

\noindent{\bf Corollary 2.26.}  {\it Assume that $M$ is projective
in $\sigma[M]$,  and let  $X\in \sigma[M]$. Then for each
submodule $Y$ of $X$, either $\sqrt[M]{Y}=X$ or $\sqrt[M]{Y}$ is
an $M$-semiprime submodule of $X$.}

\noindent {\bf Proof.}  By Theorem  2.17 and Corollary 2.25  is
clear.  $\Box$

\noindent{\bf Definition 2.27.}  Let $M$ be an $R$-module. For any
module $X$, we define ${\rm rad}_M(X)=\sqrt[M]{(0)}$. This is
called {\it $M$-Baer-McCoy radical}  or {\it $M$-prime radical} of
$X$. Thus if $X$ has an $M$-prime submodule, then ${\rm rad}_M(X)$
is equal to the intersection of all the $M$-prime submodules in
$X$ but, if $X$ has no $M$-prime submodule,
 then ${\rm rad}_M(X)=X$.

The following two propositions have  been established in
\cite{Beh1} for prime radical of modules. Now by the same method
as \cite{Beh1}, we extend these facts to $M$-prime radical of
modules.

\noindent{\bf Proposition 2.28.} {\it Let $X$ be an  $R$-module
and $Y\leq X$. Then $rad_M(Y) \subseteq rad_M(X)$.}

 \noindent
\noindent{\bf Proof.} Let $P$ be any $M$-prime submodule of $X$.
If $Y \subseteq P$, then ${\rm rad}_M(Y) \subseteq P$. If $Y
\nsubseteq P$, then it is easy to check that $Y \cap P$ is an
$M$-prime submodule of $Y$, and hence ${\rm rad}_M(Y) \subseteq (Y
\cap P) \subseteq P$. Thus in any case, ${\rm rad}_M(Y) \subseteq
P$. It follows that ${\rm rad}_M(Y) \subseteq {\rm rad}_M(X)$.
$\Box$

\noindent{\bf Lemma 2.29.} {\it Assume that $M$ is  projective
 in $\sigma[M]$, and let  $X$ be an $R$-module in
$\sigma[M]$ such that $X=\bigoplus_{\lambda\in\Lambda}X_\lambda$
is a direct sum of submodules $X_\lambda$ $(\lambda\in\Lambda)$.
Then for every submodule $N\subseteq M$, we have }{$$N\cdot
X=\bigoplus_{\lambda\in\Lambda} N\cdot {X_\lambda}$$}
 \noindent
{\bf Proof.}  Since for every $\lambda\in\Lambda$,
$X_\lambda\subseteq X$,   $N\cdot X_\lambda\subseteq N\cdot X$ for
every $\lambda\in\Lambda$. It follows that $\bigoplus_\Lambda
N\cdot {X_\lambda}\subseteq N\cdot X$. On the other hand, since
$M$ is  projective in $\sigma[M]$, so $N\cdot X=\sum_{f\in {\rm
Hom}_R(M,X)}f(N)$ and for every $\lambda\in\Lambda$, $N\cdot
{X_\lambda}=\sum_{f\in {\rm Hom}_R(M,X_\lambda)}f(N)$ by Lemma 2.9
(a). Now let $x\in N\cdot X$. Thus $x=\sum_{i=1}^{t} f_i(n_i)$
where $t\in\Bbb{N}$, $n_i\in N$ and $f_i\in {\rm Hom}_R(M,X)$.
Since $f_i(n_i)\in X$, so for every $1\leq i\leq t$,
$f_i(n_i)=\{x^{(i)}_\lambda\}_{\Lambda}$, where
$x^{(i)}_\lambda\in X_\lambda$. Thus $x=\{x^{(1)}_\lambda + ... +
x^{(t)}_\lambda\}_\Lambda=\{{\pi}_\lambda f_1(n_1) + ... +
{\pi}_\lambda f_t(n_t)\}_\Lambda$, where
${\pi}_\lambda:X\longrightarrow X_\lambda$ is the  canonical
projection for every $\lambda\in\Lambda$. It is clear that by
Lemma 2.9, $\sum_{i=1}^{t}{\pi}_{\lambda}f_i(n_i)\in N\cdot
X_{\lambda}$ for every $\lambda\in\Lambda$. Thus
$x\in\bigoplus_\Lambda N\cdot {X_\lambda}$. $\Box$

\noindent{\bf Proposition 2.30.} {\it Assume that $M$ is
projective in $\sigma[M]$,  and let $X$ be an $R$-module in
$\sigma[M]$ such that $X =\bigoplus_{\lambda\in\Lambda}X_\lambda$
is a direct sum of submodules $X_\lambda$ $(\lambda \in \Lambda)$.
Then}
$$rad_M(X)=\bigoplus_{\lambda\in\Lambda} rad_M(X_\lambda)$$
\noindent {\bf Proof.}   By  Proposition 2.28, ${\rm
rad}_M(X_\lambda)\subseteq {\rm rad}_M(X)$ for all $\lambda \in
\Lambda$. Thus  $\bigoplus_\Lambda {\rm rad}_M(X_\lambda)\subseteq
{\rm rad}_M(X)$. Now let $x \notin \bigoplus_\Lambda {\rm
rad}_M(X_\lambda)$, for some $x \in X$. Then there exists $\mu \in
\Lambda$ such that $\pi_\mu(x) \notin {\rm rad}_M(X_\mu)$, where $
\pi_\mu: X \rightarrow X_\mu$ denotes the canonical projection.
Thus there exists an $M$-prime submodule $Y_\mu$ of $X_\mu$ such
that $\pi_\mu(x) \notin Y_\mu$. Let $Z = Y_\mu \bigoplus
(\bigoplus_{\lambda\neq \mu} X_\lambda)$. It is easy to check by
Lemma 2.29 that $Z$ is an $M$-prime submodule of $X$ and $x \notin
Z$. Thus $x \notin {\rm rad}_M(X)$. It follows that ${\rm
rad}_M(X)\subseteq \bigoplus_\Lambda {\rm rad}_M(X_\lambda)$.
$\Box$

 \section{ $M$-Baer's lower nilradical of modules}

We recall the definition of the nilpotent element in a module. An
element $x$ of an $R$-module $X$ is called {\it nilpotent} if
$x=\sum_{i=1}^{r}a_{i}x_{i}$ for some $a_i \in R$, $x_i \in X$ and
$r \in\Bbb{N}$, such that ${a_i}^k{x_i}=0(1\leq i \leq r)$ for
some $k \in\Bbb{N}$ and $x$ is called {\it strongly nilpotent} if
$x=\sum_{i=1}^{r}a_{i}x_{i}$, for some $a_i \in R$, $x_i \in X$
and $r \in\Bbb{N}$, such that for every i $(1\leq i \leq r)$ and
every sequence $a_{i1},a_{i2},a_{i3}, ...$ where $a_{i1}=a_i$ and
$a_{in+1} \in a_{in}Ra_{in}(\forall n)$, we have $a_{ik}Rx_i=0$
for some $k\in\Bbb{N}$ (see [4]). It is clear that every strongly
nilpotent element of a module $X$ is a nilpotent element but the
converse is not true (see the example 2.3 [4]). In case that $R$
is commutative ring, nilpotent and strongly nilpotent are equal.

This notion has been generalized to modules over a projective module
 $M$ in $\sigma[M]$.

\noindent{\bf Definition 3.1.} Assume that $M$ is  projective  in
$\sigma[M]$,  and  let $X$  be an $R$-module in $\sigma[M]$. Then
an element $x\in X$ is called {\it $M$-nilpotent} if
$x=\sum_{i=1}^{n}r_if_i(m_i)$ for some $r_i \in R$, $m_i \in M$,
$n\in\Bbb{N}$ and $f_i \in {\rm Hom}_R(M,Rx_i)$,  where $x_i \in
X$ such that ${r_i}^kf_i(m_i)=0 (1\leq i \leq n)$ for some $k
\in\Bbb{N}$. Also, an element $x\in X$  is called {\it strongly
$M$-nilpotent}  if $x=\sum_{i=1}^{n}r_if_i(m_i)$ for some $r_i \in
R$, $m_i \in M$, $n\in\Bbb{N}$ and $f_i \in {\rm Hom}_R(M,Rx_i)$,
where $x_i \in X$ such that for every i$(1\leq i \leq n)$ and
every sequence $r_{i1},r_{i2},r_{i3},...$, where $r_{i1}=r_i$ and
$r_{it+1}\in r_{it}Rr_{it}$ $(\forall t)$, we have
$r_{ik}Rf_i(m_i)=0$ for some $k \in\Bbb{N}$.

\noindent{\bf Proposition 3.2.} {\it Let $X$ be an $R$-module.
Then an element $x \in X$ is strongly nilpotent if and only if $x$
is strongly $R$-nilpotent.}

\noindent{\bf Proof}. $(\Rightarrow )$. Suppose that  $x \in X$ is
strongly nilpotent. Then  $x=\sum_{i=1}^{n}r_ix_i$ for some $r_i
\in R$, $x_i \in X$, $n\in\Bbb{N}$ such that for every $i$ $(1\leq
i \leq n)$ and for every sequence $r_{i1},r_{i2},r_{i3},...$,
where $r_{i1}=r_i$ and $r_{it+1}\in r_{it}Rr_{it}$ $(\forall t)$,
we have $r_{ik}Rx_i=0$ for some $k\in\Bbb{N}$. Now consider
$f_i:R\rightarrow Rx_i$ such that $f_i(r)=rx_i$. Then $f_i(1)=x_i$
and it follows that
$x=\sum_{i=1}^{n}r_ix_i=\sum_{i=1}^{n}r_if_i(1)$. Since
$r_{ik}Rx_i=0$ $(1\leq i \leq n)$ for some $k \in\Bbb{N}$,we
conclude that  $r_{ik}Rf_i(1)=0$ $(1\leq i \leq n)$ for some $k
\in\Bbb{N}$,  i.e.,  $x$
is an strongly $R$-nilpotent element of $X$.\\
$( \Leftarrow )$. Assume that  $x \in X$ is strongly
$R$-nilpotent. Thus  $x=\sum_{i=1}^{n}r_if_i(a_i)$ for some $r_i,
a_i \in R$, $n \in\Bbb{N}$ and $f_i\in {\rm Hom}_R(R,Rx_i)$, where
$x_i\in X$ such that for every $i$ $(1\leq i \leq n)$ and for
every sequence $r_{i1},r_{i2},r_{i3},...$,  where $r_{i1}=r_i$ and
$r_{it+1}\in r_{it}Rr_{it}$ $(\forall t)$, we have
$r_{ik}Rf_i(a_i)=0$ for some $k \in\Bbb{N}$. Since $f_i(a_i)\in
Rx_i\subseteq X$, we conclude  that $x$ is a strongly nilpotent
 element of $X$.   $\Box$

\noindent{\bf Proposition 3.3.} {\it Let $X$ be an $R$-module.
Then an element $x \in X$ is nilpotent if and only if $x$ is
$R$-nilpotent.}

\noindent{\bf Proof}. $(\Rightarrow )$. Assume that  $x \in X$ is
nilpotent. Thus  $x=\sum_{i=1}^{n}r_ix_i$ for some $r_i \in R$,
$x_i \in X$, $n\in\Bbb{N}$ such that ${r_i}^k{x_i}=0(1\leq i \leq
n)$ for some $k \in\Bbb{N}$ . Now consider $f_i:R\rightarrow Rx_i$
such that $f_i(r)=rx_i$, so $f_i(1)=x_i$. It follows that
$x=\sum_{i=1}^{n}r_ix_i=\sum_{i=1}^{n}r_if_i(1)$. Since
${r_i}^kx_i=0$ $(1\leq i \leq n)$ for some $k \in\Bbb{N}$, so
${r_i}^kf_i(1)=0$ $(1\leq i \leq n)$ for some $k \in\Bbb{N}$,
i.e., $x$ is an $R$-nilpotent element of $X$.\\
$( \Leftarrow )$. Assume that  $x \in X$ is an  $R$-nilpotent
element. Thus  $x=\sum_{i=1}^{n}r_if_i(a_i)$ for some $r_i, a_i
\in R$, $n \in\Bbb{N}$ and $f_i\in {\rm Hom}_R(R,Rx_i)$, where
$x_i\in X$ such that ${r_i}^kf_i(a_i)=0$$(1\leq i \leq n)$ for
some $k \in\Bbb{N}$. Since $f_i(a_i)\in Rx_i\subseteq X$, we
conclude that  $x$ is a nilpotent element of $X$.$~\Box$

\noindent{\bf Proposition 3.4.} {\it Assume that $R$ is a
commutative ring,   $M$ is projective in $\sigma[M]$ and $X\in
\sigma[M]$. Then an element $x\in X$ is $M$-nilpotent if and only
if $x$ is strongly $M$-nilpotent.}

\noindent{\bf Proof}. $(\Rightarrow )$. Assume that  $x\in X$ is
$M$-nilpotent. Thus  $x=\sum_{i=1}^{n}r_if_i(m_i)$ for some $r_i
\in R$, $m_i \in M$, $n\in\Bbb{N}$ and $f_i \in {\rm
Hom}_R(M,Rx_i)$,  where $x_i \in X$ such that ${r_i}^kf_i(m_i)=0$
$(1\leq i \leq n)$ for some $k \in\Bbb{N}$. Consider sequence
$r_{i1},r_{i2},r_{i3},...$, where $r_{i1}=r_i$ and $r_{it+1}\in
r_{it}Rr_{it}$ for every $1\leq i \leq n$ and $(\forall t)$. Thus
there exists an element $r_{ik}={r_{i1}}^k{r^{\prime}}$ (where
$r^{\prime}\in R$) such that
$r_{ik}Rf_i(m_i)={r_{i1}}^k{r^{\prime}}Rf_i(m_i)=0$ (since $R$ is
commutative  and ${r_{i1}}^kf_i(m_i)=0$). Thus $x\in X$ is a strongly $M$-nilpotent element.\\
$(\Leftarrow )$. Suppose that  $x\in X$ is a strongly
$M$-nilpotent element. Thus  $x=\sum_{i=1}^{n}r_if_i(m_i)$ for
some $r_i \in R$, $m_i \in M$, $n\in\Bbb{N}$ and $f_i \in {\rm
Hom}_R(M,Rx_i)$,  where $x_i \in X$  such that for every $i$
$(1\leq i \leq n)$ and for every sequence
$r_{i1},r_{i2},r_{i3},...$, where $r_{i1}=r_i$ and $r_{it+1}\in
r_{it}Rr_{it}$ $(\forall t)$, we have $r_{ik}Rf_i(m_i)=0$ for some
$k \in\Bbb{N}$. Consider sequence $r_{i1},r_{i2},r_{i3},...$,
where $r_{i1}=r_i$ and $r_{i2}={r_{i1}}^2=r_{i1}1r_{i1}\in
r_{i1}Rr_{i1}$, $r_{i3}={r_{i1}}^4=r_{i1}1r_{i1}1r_{i1}1r_{i1}\in
r_{i2}Rr_{i2}$, ... . By assumption, we have $r_{ik}Rf_i(m_i)=0$
for some $k \in\Bbb{N}$. Since $r_{ik}={r_{i1}}^{k^\prime}$ for
some ${k^\prime}\in\Bbb{N}$, so
${r_{i1}}^{k^\prime}Rf_i(m_i)=r_{ik}Rf_i(m_i)=0$. Now for $r=1$,
we have ${r_{i1}}^{k^\prime}1f_i(m_i)=0$. Thus $x$ is an
$M$-nilpotent element. $\Box$

We recall the definition of  Baer's lower nilradical in a module.
For any module $X$, ${\rm Nil}_*(_RX)$ is the set of all strongly
nilpotent elements of $X$. In case that $R$ is a commutative ring,
${\rm Nil}_*(_RX)$ is the set of all nilpotent elements of $X$.

\noindent{\bf Definition 3.5.} Assume that  $M$ is projective
 in $\sigma[M]$. For any module $X$ in $\sigma[M]$, we
define $M$-$Nil_*(_RX)$ to be the set of all strongly
$M$-nilpotent elements of $X$. This is called {\it $M$-Baer's
lower nilradical} of $X$.

\noindent{\bf Proposition 3.6.} {\it Assume that  $M$ is
projective  in $\sigma[M]$. Then for any module X in $\sigma[M]$ }
{$$Nil_*(M).X \subseteq M{\rm -}Nil_*(_RX) \subseteq rad_M(X)$$}

\noindent{\bf Proof}. Since $M$ is  projective in $\sigma[M]$, by
Lemma 2.9 (a), $${\rm Nil}_*(M).X=\sum_{f\in {\rm
Hom}_R(M,X)}f({\rm Nil}_*(M)).$$ Now let $x\in {\rm Nil}_*(M).X$.
Thus $x=\sum_{i=1}^{s}f_i(m_i)$ for some $m_i\in {\rm Nil}_*(M)$,
$s\in\Bbb(N)$ and $f_i\in {\rm Hom}_R(M,X)$. Since $m_i\in {\rm
Nil}_*(M)$, so $m_i=\sum_{j=1}^{t}r_{i_j}n_{i_j}$ for some
$r_{i_j}\in R$, $n_{i_j}\in M$, $t\in\Bbb{N}$ such that for every
$j$ $(1\leq j \leq t)$ and for every sequence
$r_{i_{j1}},r_{i_{j2}},r_{i_{j3}},...$,  where
$r_{i_{j1}}=r_{i_j}$ and $r_{i_{ju+1}}\in r_{i_{ju}}Rr_{i_{ju}}$
$(\forall u)$, we have $r_{i_{j{k_i}}}Rn_{i_j}=0$ for some $k_i
\in\Bbb{N}$. Thus
$x=\sum_{i=1}^{s}f_i(m_i)=\sum_{i=1}^{s}f_i(\sum_{j=1}^{t}r_{i_j}n_{i_j})=\sum_{i=1}^{s}\sum_{j=1}^{t}r_{i_j}f_i(n_{i_j})$.
Since $r_{i_{j{k_i}}}Rn_{i_j}=0$, we conclude  that
$0=f_i(r_{i_{j{k_i}}}Rn_{i_j})=r_{i_{j{k_i}}}Rf_i(n_{i_j})$
 for some $k_i\in\Bbb{N}$, where $(1\leq
i \leq s)$ and  $(1\leq j \leq t)$. Thus $x\in M$-${\rm
Nil}_*(_RX)$.\\
 \indent Let $x\in M$-${\rm Nil}_*(_RX)$ and $x\notin {\rm
rad}_M(X)={\sqrt[M]{(0)}}$. So $x=\sum_{i=1}^{n}a_if_i(m_i)$ for
some $a_i\in R$, $m_i\in M$, $n\in\Bbb{N}$ and $f_i\in {\rm
Hom}_R(M,Rx_i)$ such that for every i$(1\leq i \leq n)$ and for
every sequence $a_{i1},a_{i2},a_{i3},...$,  where $a_{i1}=a_i$ and
$a_{iu+1}\in a_{iu}Ra_{iu}$ $(\forall u)$, we have
$a_{ik}Rf_i(m_i)=0$ for some $k \in\Bbb{N}$. Without loss of
generality, we can assume that $a_1f_1(m_1)\notin {\rm rad}_M(X)$.
Thus there exists an $M$-$m$-system $S$ such that $a_1f_1(m_1)\in
S$ and $0\notin S$. On the other hand $a_1f_1(m_1)\in
Ra_1(Rm_1).(Rx_1)$. Thus $Ra_1(Rm_1).(Rx_1)~\cap S\neq\emptyset$
and hence $Ra_1(Rm_1).X~\cap S\neq\emptyset$. Therefore,  if we
put $N=Ra_1(Rm_1)$, $Y=(0)$ and $Z=Ra_1(Rm_1).(Rx_1)$, then
$(Ra_1(Rm_1))^2.(Rx_1)~\cap S\neq\emptyset$  by Proposition 2.14.
Since $M$ is projective in $\sigma[M]$, by Lemma 2.9 (a)
and Lemma 2.23,  we conclude that \vspace{2mm}\\
$~~~~~~~~~~~~~~(Ra_1(Rm_1))^2.(Rx_1)=(Ra_1(Rm_1).Ra_1(Rm_1)).(Rx_1)\\~~~~~~~~~~~~~~~~~~~~~~~~~~~~~~~~~~~~~~~~=(Ra_1(Rm_1)).(Ra_1(Rm_1).(Rx_1))\\
~~~~~~~~~~~~~~~~~~~~~~~~~~~~~~~~~~~~~~~~=\sum_{f\in
{\rm Hom}_R(M,Ra_1(Rm_1).(Rx_1))}f(Ra_1(Rm_1)).$\vspace{2mm}\\
\indent Assume that $s_1=1$, $a_{11}=a_1$ and
$a_1f_1(t_1a_1s_2m_1)\in (Ra_1(Rm_1))^2.(Rx_1)~\cap S$, where
$s_2,t_1\in R$. Since $a_1f_1(t_1a_1s_2m_1)=s_2a_1t_1a_1f_1(m_1)$
and $a_{12}=a_1t_1a_1$, so   $s_2a_{12}f_1(m_1)\in
Ra_{12}(Rm_1).(Rx_1)\cap S$. It follows that
$Ra_{12}(Rm_1).(Rx_1)~\cap S\neq\emptyset$ and so
$$(Ra_{12}(Rm_1))^2.(Rx_1)~\cap S\neq\emptyset.$$ Thus there exists
$s_3a_{13}f_1(m_1)\in (Ra_{12}(Rm_1))^2.(Rx_1)\cap S$, where
$s_3\in R$, and $a_{13}:=a_{12}t_2s_2a_{12}$ for some $t_2\in R$.
We can repeat this argument to get sequences $\{s_u\}_{u\in N}$
and $\{a_{1u}\}_{u\in\Bbb{N}}$ in $R$, where $a_{11}=a_1$ and
$a_{1u+1}\in a_{1u}Ra_{1u}$ $(\forall u)$, such that
$s_ua_{1u}f_1(m_1)\in S$ for all $u\geq 1$. Now by our hypothesis
$a_{1k}Rf_1(m_1)=0$ for some $k\in\Bbb{N}$, and so
$s_ka_{1k}f_1(m_1)=0\in S$, a contradiction.        $\Box$

In case $M=R$, by Proposition 3.6,  ${\rm Nil}_*(R).X \subseteq
R$-${\rm Nil}_*(_RX) \subseteq {\rm rad}_R(X)$. Since by
Proposition 3.2,  $R$-${\rm Nil}_*(_RX)$ is the set of all
strongly $R$-nilpotent elements of $X$, so we have $R$-${\rm
Nil}_*(_RX)={\rm Nil}_*(_RX)$ (see also, \cite[Lemma 3.2]{Beh1}).

\noindent{\bf Corollary 3.7.} {\it Assume that  $M$ is projective
 in $\sigma[M]$.  Then $$Nil_*(M)=Nil_*(M).M=M-Nil_*(M).$$}
\noindent{\bf Proof}. By Proposition 3.6, ${\rm Nil}_*(M).M
\subseteq M$-${\rm Nil}_*(M)$. Also, we have ${\rm
Nil}_*(M).M=\sum_{f\in {\rm Hom}_R(M,M)}f({\rm Nil}_*(M))$, by
Lemma 2.9 (a). Since $1_M\in {\rm Hom}_R(M,M)$, so ${\rm
Nil}_*(M)\subseteq {\rm Nil}_*(M).M$. On the other hand, if $x\in
M$-${\rm Nil}_*(M)$, then  $x=\sum_{i=1}^{n}r_if_i(m_i)$ for some
$r_i\in R$, $m_i\in M$, $n\in\Bbb{N}$ and $f_i\in {\rm
Hom}_R(M,Rx_i)$, where $x_i\in M$ such that for every $i$ $(1\leq
i \leq n)$ and for  every sequence $r_{i1},r_{i2},r_{i3},...$,
where $r_{i1}=r_i$ and $r_{it+1}\in r_{it}Rr_{it}$ $(\forall t)$,
we have $r_{ik}Rf_i(m_i)=0$ for some $k \in\Bbb{N}$. Since
$f_i(m_i)\in Rx_i\subseteq M$, it follows that $x$ is a strongly
nilpotent element of $M$. So $x\in {\rm Nil}_*(M)$. It follows
that $M$-${\rm Nil}_*(M)\subseteq {\rm Nil}_*(M)$ and ${\rm
Nil}_*(M)\subseteq {\rm Nil}_*(M).M\subseteq M$-${\rm
Nil}_*(M)\subseteq {\rm Nil}_*(M)$.
 Thus  ${\rm Nil}_*(M)={\rm Nil}_*(M).M=M$-${\rm Nil}_*(M)$.   $\Box$

\noindent{\bf Corollary 3.8.} {\it Assume that  $M$ is projective
 in $\sigma[M]$. Then $rad_R(M)\subseteq rad_M(M)$.}

\noindent{\bf Proof}. By Proposition 3.6, we have $M$-${\rm
Nil}_*(M) \subseteq {\rm rad}_M(M)$. On the other hand ${\rm
Nil}_*(M)=M$-${\rm Nil}_*(M)$ by Corollary 3.7. Thus  ${\rm
Nil}_*(M)\subseteq {\rm rad}_M(M)$. Since $M$ is  projective in
$\sigma[M]$,  ${\rm rad}_R(M)={\rm Nil}_*(M)$ by \cite[Theorem
3.8]{Beh1}. Thus ${\rm rad}_R(M)={\rm Nil}_*(M)\subseteq {\rm
rad}_M(M)$.$~\Box$

\noindent{\bf Proposition 3.9.} {\it Assume that  $M$ is
projective    in $\sigma[M]$. If $X\in\sigma[M]$ such that
$rad_M(X)=M$-$Nil_*(X)$, then $rad_M(Y)=M$-$Nil_*(Y)$ for any
direct summand $Y$ of $X$.}

\noindent{\bf Proof}. Suppose that $X=Y\oplus Z$, where  $Z$, $Y$
are submodules  of $X$. By Proposition 3.6, $M$-${\rm
Nil}_*(Y)\subseteq {\rm rad}_M(Y)$. Let $x\in {\rm rad}_M(Y)$. By
Proposition 2.28, $x\in {\rm rad}_M(X)$. By hypothesis $x\in
M$-${\rm Nil}_*(X)$. Thus $x=\sum_{i=1}^{n}r_if_i(m_i)$ for some
$r_i \in R$, $m_i \in M$, $n\in\Bbb{N}$ and $f_i \in {\rm
Hom}_R(M,Rx_i)$,  where $x_i \in X$ such that for every $i$
$(1\leq i \leq n)$ and for every sequence
$r_{i1},r_{i2},r_{i3},...$, where $r_{i1}=r_i$ and $r_{it+1}\in
r_{it}Rr_{it}$ $(\forall t)$, we have $r_{ik}Rf_i(m_i)=0$ for some
$k \in\Bbb{N}$. Since $x_i\in X$, there exist elements $y_i\in Y
$, $z_i\in Z$ such that $x_i=y_i+z_i$ for each $i$ $(1\leq i \leq
n)$. On the other hand, $f_i(m_i)\in Rx_i$ for each $i$, and hence
 $f_i(m_i)=a_i(y_i+z_i)$ for some $a_i\in R$ $(1\leq i
\leq n)$. It is clear that $x=r_1a_1y_1+r_2a_2y_2+...+r_na_ny_n$,
and $r_{ik}Ra_iy_i=0$, for some $k\in\Bbb{N}$,$(1\leq i \leq n)$.
Now for each $i$ $(1\leq i \leq n)$, we consider
$g_i:M\stackrel{f_i}\longrightarrow Rx_i\subseteq
X\stackrel{{\pi}_i}\longrightarrow Ry_i\subseteq Y$, where
${\pi}_i$ is the natural projection map such that
$g_i(m_i)={\pi}_if_i(m_i)={\pi}_i(a_i(y_i+z_i))=a_iy_i$. Thus
$x=r_1a_1y_1+r_2a_2y_2+...+r_na_ny_n=\sum_{i=1}^{n}r_ig_i(m_i)$,
where $g_i\in {\rm Hom}_R(M,Ry_i)$ and
$r_{ik}Ra_iy_i=r_{ik}Rg_i(m_i)=0$. It follows that $x\in M$-${\rm
Nil}_*(Y)$. Thus  ${\rm rad}_M(Y)=M$-${\rm Nil}_*(Y)$. $\Box$

\section{$M$-injective modules and prime $M$-ideals  }

The module $_RX$ is said to be {\it $M$-generated} if there exists
an $R$-epimorphism from a direct sum of copies of $M$ onto $X$.
Equivalently, for each nonzero $R$-homomorphism $f: X\rightarrow
Y$ there exists an $R$-homomorphism $g: M \rightarrow X$ with
$fg\neq 0$. The {\it trace} of $M$ in $X$ is defined to be
$$tr^M(X)=\sum_{f\in {\rm Hom}_R(M,X)} f(M)$$
 and thus $X$ is $M$-generated if and only if $tr^M(X)=X$.

We recall the definition  of prime $M$-ideal.  The proper
$M$-ideal $P$ is said to be a prime $M$-ideal if there exists an
$M$-prime module $_RX$ such that $P = {\rm Ann}_M(X)$.

\noindent{\bf Proposition 4.1.} {\it Let $M$ be an $R$-module with
$Hom_R(M,X)\neq 0$ for every $X\in\sigma[M]$ and $P$ be a proper
$M$-ideal. Then $P$ is a prime $M$-ideal if and only if $P$ is a
Beachy-prime $M$-ideal.}

 \noindent {\bf Proof.} Assume that $P$ is  a prime $M$-ideal. Thus  there
 exists $M$-prime module $X$ such that $P={\rm Ann}_M(X)$. Since $P\neq
 M$,  ${\rm Hom}_R(M,X)\neq 0$. Thus by  Proposition 2.7,
  $X$ is a Beachy-$M$-prime module. Thus $P$ is a Beachy-prime
 $M$-ideal.\\
 \indent Conversely, let $P$ be a Beachy-prime $M$-ideal. Thus  there
 exists a Beachy-$M$-prime module $X$ in $\sigma[M]$ such that $P={\rm Ann}_M(X)$.
 Since ${\rm Hom}_R(M,X)\neq0$, so $X\neq (0)$. Now
 assume that $Y$ is  a nonzero submodule of $X$. So $Y\in\sigma[M]$ and
 ${\rm Hom}_R(M,Y)\neq 0$ by assumption. Therefore,  ${\rm Ann}_M(X)={\rm Ann}_M(Y)$ by
 the definition of Beachy-$M$-prime module. Thus by Proposition 2.4, $X$ is an $M$-prime
 module and hence $P$ is a prime $M$-ideal.  $\Box$

 The module $_RX$ in $\sigma[M]$ is said to be {\it finitely $M$-generated} if
there exists an epimorphism $f : M^n\rightarrow X$, for some
positive integer $n$. It is said to be {\it finitely
$M$-annihilated}  if there exists a monomorphism $g: M/{\rm
Ann}_M(X)\rightarrow X^m$, for some positive integer $m$. Also,
the module $_RM$ is said to {\it satisfy condition $H$} if every
finitely $M$-generated module is finitely $M$-annihilated. Note
that if $M=R$ and $R$ is a fully bounded Noetherian ring, then $M$
satisfies condition $H$. The same is true if $M$ is an Artinian
module, since then $M/K$ has the finite intersection property.

In [1, Theorem 6.7], it is shown that if  $M$ is  a Noetherian
module such that $M$ satisfies condition $H$ and ${\rm
Hom}_R(M,X)\neq 0$ for all modules $X$ in $\sigma[M]$, then there
is a one-to-one correspondence between isomorphism classes of
indecomposable $M$-injective modules in $\sigma[M]$ and
Beachy-prime $M$-ideals. Next, in the main result  of this
section, we show   this fact is also true for a Noetherian module
with  condition $H$ and the assumption ${\rm Hom}_R(M,X)\neq 0$
for all modules $X$ in $\sigma[M]$
 via prime $M$-ideals.

\noindent{\bf Corollary 4.2.} {\it Let $M$ be a Noetherian
$R$-module. If $M$ satisfes condition $H$ and $Hom_R(M,X)\neq 0$
for all modules $X$ in $\sigma[M]$, then there is a one-to-one
correspondence between isomorphism classes of indecomposable
$M$-injective modules in $\sigma[M]$ and prime $M$-ideals.}

\noindent {\bf Proof.} By  \cite[Theorem 6.7]{Beachy} and
Proposition 4.1, is clear.

\section{Prime $M$-ideals and $M$-prime radical of Artinian modules  }

Let $M$ be an $R$-module. Recall that a proper submodule $P$ of
$M$ is {\it virtually maximal} if the factor module $M/P$ is a
homogeneous semisimple $R$-module,  i.e., $M/P$ is a direct sum of
isomorphic simple modules. Clearly, every virtually maximal
submodule of $M$ is prime. Also, every maximal submodule of $M$ is
virtually maximal and for $M=R$ and $R$ commutative, this is
equivalent to the notion of maximal ideal in $R$.

We recall that ${\rm Soc}(M)$ is sum of all minimal submodules of
$M$. If $M$ has no minimal submodule, then ${\rm Soc}(M)=(0)$.

\noindent{\bf Proposition 5.1.} {\it Let $M$ be an Artinian
R-module. If $M$ is an $M$-prime module, then $M$ is a homogeneous
semisimple module.}

\noindent {\bf Proof.}   Since $M$ is an Artinian $R$-module,
${\rm Soc}(M)\neq (0)$. Hence there exist simple submodule $Rm$ of
$M$ where $0\neq m\in M$. Since $M$ is an $M$-prime module, ${\rm
Ann}_M(Rm)={\rm Ann}_M(M)=(0)$ by Proposition 2.4. Thus $(0)={\rm
Ann}_M(Rm)=\bigcap_{f\in {\rm Hom}_R(M,Rm)}{\rm {\rm ker}}(f)$.
Since $Rm\cong M/{\rm {\rm ker}}(f)$ for every $f\in {\rm
Hom}_R(M,Rm)$,  $(0)$ is an intersection of maximal submodules and
since $M$ is Artinian, $(0)$ must be a finite intersection of
maximal submodules. It follows that $M$ is isomorphic to a finite
direct sum of copies of $Rm$. Thus  $M$ is a homogeneous
semisimple module. $\Box$

An $M$-ideal $P$ is said to be a {\it primitive $M$-ideal} if
$P={\rm Ann}_M(S)$ for a simple module $_RS$ (see \cite[Definition
3.5]{Beachy}).

\noindent{\bf Proposition 5.2.} {\it Let $P$ be a proper
$M$-ideal. If $P$ is a primitive $M$-ideal, then $P$ is a prime
$M$-ideal.}

\noindent {\bf Proof.}   If $P$ is a primitive $M$-ideal, then
$P={\rm Ann}_M(S)$ for a simple  $R$-module $S$. Since $S$ has no
any proper submodule,  $S$ is an $M$-prime module by Proposition
2.4. Thus $P$ is a prime $M$-ideal. $\Box$

\noindent{\bf Proposition 5.3.} {\it Let $M$ be an $M$-prime
module with $Soc(M)\neq (0)$. Then $(0)$ is a primitive
$M$-ideal.}

\noindent {\bf Proof.}   Since ${\rm Soc}(M)\neq (0)$,  there
exists a simple submodule $Rm$ of $M$ where $0\neq m\in M$. Since
$M$ is an $M$-prime module, so ${\rm Ann}_M(Rm)={\rm
Ann}_M(M)=(0)$. Therefore,  $(0)$ is a primitive $M$-ideal. $\Box$

\noindent{\bf Proposition 5.4.} {\it Assume that  $M$ is
projective in $\sigma[M]$. If $M$ is an Artinian $R$-module, then
every prime $M$-ideal of $M$ is virtually maximal.}

\noindent {\bf Proof.} Suppose that  $P\lneqq M$ is  a prime
$M$-ideal. Since $M$ is projective in $\sigma[M]$,  $M/P$ is an
$M$-prime module by Proposition 2.10. Since $M/P$ is also an
Artinian module,  ${\rm Soc}(M/P)\neq (0)$ and hence there exists
a simple submodule $R\bar{m}$ of $M/P$ where $0\neq\bar{m}\in
M/P$. Since $M/P$ is an $M$-prime module,  ${\rm
Ann}_M(R\bar{m})={\rm Ann}_M(M/P)=P$. On the other hand, $P={\rm
Ann}_M(R\bar{m})=\bigcap_{f\in {\rm Hom}_R(M,R\bar{m})}{\rm
ker}(f)$. Since $R\bar{m}\cong M/{\rm ker}(f)$ for every $f\in
{\rm Hom}_R(M,R\bar{m})$,  $P$ must be an intersection of maximal
submodules. Since $M/P$ is Artinian, $P$ must be a finite
intersection of maximal submodules, and so  $M/P$ is isomorphic to
a finite direct sum of copies of $R\bar{m}$. Thus  $M/P$ is a
homogeneous semisimple module,  i.e.,  $P$ is a virtually maximal
submodule of $M$.      $\Box$\

\noindent{\bf Definition 5.5.} The {\it prime radical} of the
module $M$, denoted by $P(M)$, is defined to be the intersection
of all prime $M$-ideals.

We note that each prime $M$-ideal is the annihilator  of an
$M$-prime module in $M$. It follows that $P(M)={\rm rad}_{\cal
C}(M)$, where $\cal C$ is the class of all $M$-prime left
$R$-modules. If $_RX$ is any module with a submodule $Y$ such that
$X/Y$ is an $M$-prime module, then ${\rm rad}_{\cal C}(X)\subseteq
Y$. In this case it follows from \cite[Lemma 1.8]{Beachy} that
$P(M).X\subseteq Y$.

\noindent{\bf Theorem 5.6.} {\it Assume that  $M$ is projective in
$\sigma[M]$. If $M$ is an Artinian $R$-module, then every prime
$M$-ideal of $M$ is virtually maximal and $M/P(M)$ is a Noetherian
$R$-module.}

\noindent {\bf Proof.} If $M$ does not contain any prime
$M$-ideal, then $P(M)=M$. Suppose that $M$ contains a prime
$M$-ideal. By Proposition 5.4, every prime $M$-ideal of $M$ is
virtually maximal. Let $N$ be minimal in the collection $\cal S$
of $M$-ideals of $M$ which are finite intersections of primes. If
$P$ is any prime $M$-ideal of $M$, then $P\cap N\in\cal S$ and
$P\cap N\subseteq N$. Thus  $N=P\cap N\subseteq P$ by minimality
of $N$ in $\cal S$. It follows that $N=P(M)$. On the other hand,
for each prime $M$-ideal, the factor module $M/P$ is a homogeneous
semisimple module with DCC. So $M/P$ is Noetherian. Thus $M/P$ is
Noetherian for every prime $M$-ideal $P$ of $M$. Since $P(M)$ is a
finite intersection of prime $M$-ideals,  $M/P(M)$ is  also a
Noetherian $R$-module.      $\Box$

The following theorem is a generalization of \cite[Theorem
2.11]{Beh1}.

\noindent{\bf Theorem 5.7.} {\it Assume that  $M$ is projective in
$\sigma[M]$. If $M$ be an Artinian $R$-module, then $P(M)=M$ or
there exist primitive $M$-ideals $P_1$,...,$P_n$ of $M$ such that
$P(M)=\bigcap_{i=1}^{n}P_i$.}

\noindent {\bf Proof.}   Let $P$ be a prime $M$-ideal of $M$.
Since $M$ is projective in $\sigma[M]$, so $M/P$ is an $M$-prime
module by Proposition 2.10 (ii). Since $M/P$ is an Artinian
$R$-module, ${\rm Soc}(M/P)\neq (0)$. Thus  there exists a simple
submodule $R\bar{m}$ of $M/P$ where $0\neq \bar{m}\in M/P$. Since
$M/P$ is an $M$-prime module,  ${\rm Ann}_M(R\bar{m})={\rm
Ann}_M(M/P)$. On the other hand, ${\rm Ann}_M(M/P)=P$, since $P$
is an $M$-ideal. Thus $P$ is a primitive $M$-ideal. Since $P$ is
arbitrary prime $M$-ideal, so every prime $M$-ideal of $M$ is
primitive $M$-ideal. On the other hand by Proposition 5.2, we have
that every primitive $M$-ideals is prime $M$-ideal. Thus  $P(M)$
is the intersection all of primitive $M$-ideal of $M$. Now let $N$
be minimal in the collection $\cal S$ of $M$-ideals of $M$ which
are finite intersections of primes. If $Q$ is any prime $M$-ideal
of $M$, then $Q\cap N\in\cal S$ and $Q\cap N\subseteq N$. Thus
$N=Q\cap N\subseteq Q$ by minimality of $N$ in $\cal S$. It
follows that $N=P(M)$. Thus $P(M)$ is a finite intersection of
prime $M$-ideals and it follows that $P(M)$ is a finite
intersection of primitive $M$-ideals. So there exist primitive
$M$-ideals $P_1$,...,$P_n$ of $M$ such that
$P(M)=\bigcap_{i=1}^{n}P_i$. Since $P_i$ is an $M$-ideal for every
$1\leq i\leq n$,  $P_i.M=P_i$ and so
$P(M)=\bigcap_{i=1}^{n}P_i.M=\bigcap_{i=1}^{n}P_i$.       $\Box$

\noindent{\bf Corollary 5.8.} {\it  Assume that  $M$ is projective
in $\sigma[M]$. If $M$ be an Artinian $M$-prime module, then
$P(M)=(0)$.}

\noindent {\bf Proof.} By Proposition 5.1,  $(0)$ is a primitive
$M$-ideal of $M$. It follows that $P(M)=(0)$ by Theorem 5.8.
$\Box$

Minimal $M$-prime submodules are defined in a natural way. By
Zorn's Lemma one can easily see that each $M$-prime submodule of a
module $X$ contains a minimal $M$-prime submodule of $X$. In
\cite[Theorem 5.2]{MS},   it is shown that every  Noetherian
module contain only finitely many minimal prime submodules. It is
easy to show that if $X$ is a Noetherian module, then $X$ contain
only finitely many minimal $M$-prime submodules.

We conclude this paper with the following interesting result,
which is  a generalization of \cite[Theorem 2.1]{Beh1}.

 \noindent{\bf Theorem 5.9.} {\it  Let $X$ be a Noetherian
 $R$-module. If every $M$-prime submodule of $X$ is virtually maximal,
 then $X/rad_M(X)$ is an Artinian  $R$-module.}

 \noindent{\bf Proof.} By our
hypotheses, for each $M$-prime submodule $P$ of $X$,   $X/P$ is  a
homogeneous semisimple  $R$-module. Since $X$ is a Noetherian
$R$-module, $X/P$ is also  Noetherian. This implies that $X/P$ is
an Artinian $R$-module. On the other hand ${\rm rad}_M(X)=P_1\cap
\cdots \cap P_n$ where $P_1, \cdots, P_n$ are all minimal
$M$-prime submodules of $M$. Thus $X/P_1\oplus \cdots\oplus X/P_n$
is also an Artinian  $R$-module. It follows that
 $X/{\rm rad}_M(X)$ is  an Artinian  $R$-module.  $\square$

\end{document}